\documentclass[10pt,a4paper,reqno]{amsart} 

\usepackage{graphicx,psfrag,
amsfonts,amsmath,indentfirst,floatflt}
\usepackage{graphics,color}

\usepackage[T1]{fontenc}
 \usepackage{latexsym}
\usepackage{rotating,epsfig}
\usepackage{amssymb}
\usepackage{amsfonts}

\newtheorem{Theorem}{Theorem}
\newtheorem{Lemma}[Theorem]{Lemma}

\begin{document}


 \newfont{\bbold}{msbm10}
 \newfont{\bbolds}{msbm7}

\newcommand{\epsf}[1]{\epsfbox{#1}}
\newcommand{\nss}{\mbox{\bbolds N}}
\newcommand{\zss}{\mbox{\bbolds Z}}

\newcommand{\ns}{\mbox{\bbold N}}
\newcommand{\zs}{\mbox{\bbold Z}}
\newcommand{\qs}{\mbox{\bbold Q}}
\newcommand{\rs}{\mbox{\bbold R}}
\newcommand{\cs}{\mbox{\bbold C}}

\newcommand{\bm}[1]{\mbox{\boldmath \ensuremath{#1}}}
\newcommand{\bs}[1]{\mbox{\boldmath \ensuremath{\scriptstyle #1}}}

\newcommand{\bx}{\bar x}
\newcommand{\by}{\bar y}
\newcommand{\bu}{\bar u}

\newcommand{\GL}{\mathbb{L}}
\newcommand{\GM}{\mathbb{M}}
\newcommand{\GK}{\mathbb{K}}
\newcommand{\GP}{\mathbb{P}}
\newcommand{\GD}{\mathbb{D}}
\newcommand{\GE}{\mathbb{E}}
\newcommand{\GA}{\mathbb{A}}
\newcommand{\GS}{\mathbb{S}}
\newcommand{\GC}{\mathbb{C}}
\newcommand{\GF}{\mathbb{F}}
\newcommand{\GB}{\mathbb{B}}
\newcommand{\GR}{\mathbb{R}}
\newcommand{\GQ}{\mathbb{Q}}
\newcommand{\GU}{\mathbb{U}}
\newcommand{\GW}{\mathbb{W}}

\newcommand{\Ref}[1]{(\ref{#1})}
\newcommand{\beq}{\begin{equation}}
\newcommand{\eeq}{\end{equation}}
\newcommand{\gf}{generating function}
\newcommand{\sg}{s\'erie g\'en\'eratrice}
\newcommand{\gfs}{generating functions}
\newcommand{\al}{\alpha}
\newcommand{\be}{\beta}
\newcommand{\cd}{\cdot}
\def\cqfd{\par\nopagebreak\rightline{\vrule height 3pt width 5pt depth 2pt}
\medbreak}

\title[Counting walks in the quarter plane]{Counting walks in the quarter plane}



\author[M. Bousquet-M\'elou]{Mireille Bousquet-M\'elou}

\address{CNRS, LaBRI, Universit\'e de Bordeaux, 351 cours de la
   Lib\'eration,  33405 Talence Cedex, France}

\maketitle

\begin{abstract}
We study planar walks that start from a given point $(i_0,j_0)$, take
their steps in a finite set $\frak S$, and are confined in the first
quadrant $x \ge 0, y\ge 0$. Their enumeration can be attacked in a
systematic way: the \gf \ $Q(x,y;t)$ that counts them by their length
(variable $t$) and the coordinates of their endpoint (variables $x,y$)
satisfies a linear 
functional equation encoding the step-by-step description of
walks. For instance, for the square lattice walks starting from the
origin, this equation reads 
$$ \left(xy -t(x+y+x^2y+xy^2)  \right) Q(x,y;t) = xy-xtQ(x,0;t)
-ytQ(0,y;t).$$
The central question addressed in this paper is the {\em nature\/} of
the series $Q(x,y;t)$. When is it algebraic? When is it D-finite (or 
holonomic)? Can these 
properties be derived from the functional equation itself?

Our first result is a new proof of an old theorem due to Kreweras,
according to which one of these walk models has, for mysterious
reasons, an algebraic \gf . Then, we provide a new proof of a holonomy
criterion recently proved by M.~Petkov\v sek and the author. In both
cases, we work directly from the functional equation.

\end{abstract}



\section{\bf Walks in the quarter plane}
The enumeration of lattice walks is one of the most venerable topics
in enumerative combinatorics, which has numerous applications in
probability~\cite{feller,lawler,spitzer}. These walks take their steps in a
finite subset $\frak S$ of $\zs ^d$, and might be constrained in
various ways. One can only cite a small percentage of the relevant
litterature, which dates back at least to the next-to-last century
~\cite{andre,gessel-proba,kreweras,mohanty,narayana}.
Many recent publications show that the topic is still
active~\cite{firenze,mbm-slitplane,prep-GS,gessel-zeilberger,grabiner,niederhausen1,niederhausen2}.

After the solution of many explicit problems, certain patterns have
emerged, and a more recent trend consists in developing methods that
are valid for generic sets of steps.
A special attention is being paid
to the {\em nature\/} of the \gf \ of the walks under
consideration. For instance, the \gf \ for unconstrained walks on the
line $\zs$ is rational, while the \gf \ for walks constrained to stay
in the half-line $\ns$ is always algebraic~\cite{banderier-flajolet}.
This result has often been described in terms of {\em partially directed\/}
2-dimensional walks confined in a quadrant (or {\em generalized Dyck
walks\/}~\cite{duchon,gessel-factor,lab3,lab2}), but is, essentially, of a 
1-dimensional nature.

Similar questions can be addressed for  {\em real\/} 2-dimensional
walks. Again, the \gf \ for unconstrained walks starting 
from a given 
point is clearly rational. Moreover, the argument used for
1-dimensional walks confined in $\ns$ can
be recycled to prove that the \gf \ for the  walks that
stay in the half-plane $x\ge 0$ is always algebraic. What about
doubly-restricted walks, that is, walks that are confined in the
quadrant $x\ge 0, y \ge 0$?  

\begin{figure}[ht]
\begin{center}
\includegraphics[height=40mm]{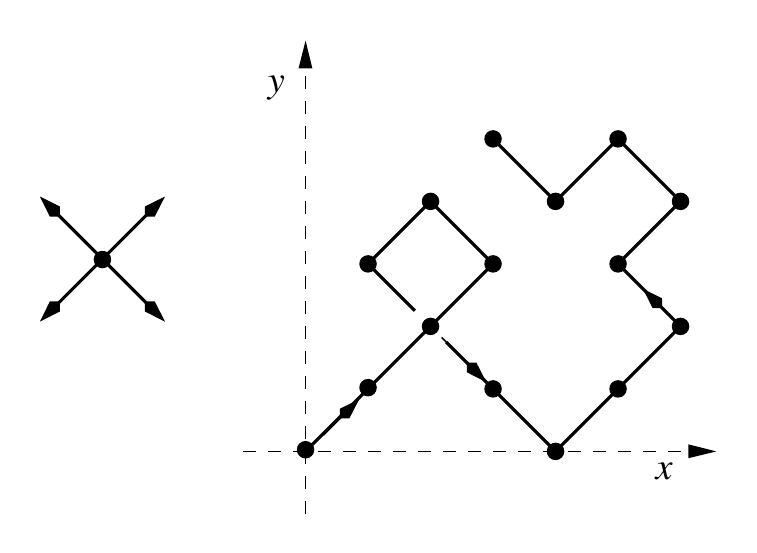}
\end{center}
\caption{A walk on the diagonal square lattice confined in the first quadrant.}
\label{diagonal}
\end{figure}

A rapid inspection of the most standard 
cases suggests that these walks might have always a D-finite \gf
\footnote{A series $F(t)$ is 
D-finite (or {\em holonomic\/}) if it satisfies a linear differential equation with 
polynomial coefficients in $t$. Any algebraic series is
D-finite.}. The simplest example is probably that of the diagonal square
lattice, where the steps are 
North-East, South-East, North-West and South-West
(Figure~\ref{diagonal}): by projecting the walks on the $x$- and
$y$-axes, we 
obtain two decoupled prefixes of Dyck paths, so that the length \gf \ for
walks that start from the origin and stay in the first quadrant is
$$\sum_{n \ge 0} {n \choose {\lfloor n/2\rfloor}}^2 t^n,$$
 a D-finite series. 
For the ordinary square lattice (with North, East, South and West steps),
the \gf \ is
$$\sum_{ m, n \ge 0} {m+n \choose m}
{m \choose {\lfloor m/2\rfloor}} 
{n \choose {\lfloor n/2\rfloor}}  t^{m+n}
=
\sum_{ n \ge 0} {n \choose {\lfloor n/2\rfloor}} 
{{n+1} \choose {\lceil n/2\rceil}} t^{n},$$
another D-finite series. The first expression comes from the fact that
these walks are  shuffles of two prefixes of Dyck walks, and the
Chu-Vandermonde identity transforms it into the second simpler
expression~\cite{kratt}.

In both cases, the number of $n$-step walks  grows
asymptotically like $4^n/n$, which prevents the \gf \
from being algebraic (see~\cite{flajolet} for the possible asymptotic
behaviours of coefficients of algebraic series).

The two above results can be refined by taking into account the
coordinates  of the endpoint: if $a_{i,j}(n)$ denotes the number of
$n$-step walks of length $n$ ending at $(i,j)$, then we have, for the
diagonal square lattice:
$$
\sum_{i,j,n \ge 0} a_{i,j}(n) x^i y^j t^n = 
\sum_{i,j,n \ge 0} \frac{(i+1)(j+1)}{(n+1)^2} 
{{n+1} \choose {\frac {n-i} 2}}
 {{n+1} \choose {\frac {n-j} 2}}
x^i y^j t^{n},$$
where the binomial coefficient 
$ {n \choose {( {n-i})/ 2}}$
is zero unless $0\le i \le n$ and $i\equiv n$ mod $2$.
Similarly, for the ordinary square lattice,
\begin{equation}
\sum_{i,j,n \ge 0} a_{i,j}(n) x^i y^j t^n = 
\sum_{i,j,n \ge 0} \frac{(i+1)(j+1)}{(n+1)(n+2)} 
{{n+2} \choose {\frac {n+i-j+2} 2}}
 {{n+2} \choose {\frac {n-i-j} 2}}
x^i y^j t^{n}.
\label{carre-general}
\end{equation}
%
%
%
These two series can be seen to be D-finite in their three variables.

\medskip

This holonomy, however, is not the rule: as proved
in~\cite{mbm-marko2},  walks that start from $(1,1)$, take their steps
in ${\frak S}= 
\{(2,-1), (-1,2)\}$ and always stay in the first quadrant have a
non-D-finite length \gf . The same holds for the subclass of walks
ending on the $x$-axis.  These walks are sometimes called {\em
knight's walks\/}.

\medskip

\begin{figure}[ht]
\begin{center}
\includegraphics[height=40mm]{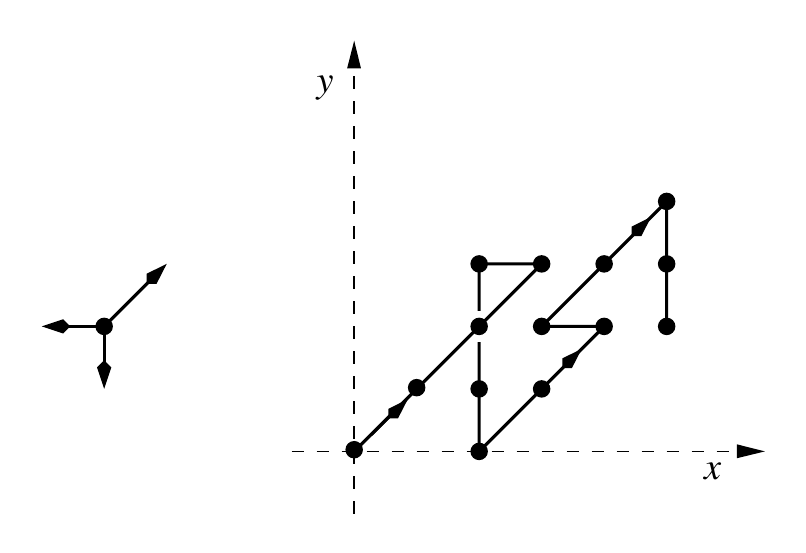}
\end{center}
\caption{Kreweras' walks in a quadrant.}
\label{figure-kreweras}
\end{figure}

At the other end of the hierarchy, another walk model displays a
mysteriously simple algebraic \gf : when the starting point is
$(0,0)$, and the allowed steps South, West and North-East
(Figure~\ref{figure-kreweras}), the number  
of walks of length $3n+2i$ ending at the point $(i,0)$ is
\begin{equation} 
\frac{4 ^n (2i+1)}{(n+i+1)(2n+2i+1)} 
{{2i} \choose i}
{{3n+2i} \choose n} .
\label{kreweras-simple}\end{equation}
This result was first proved by Kreweras in 1965~\cite[Chap.~3]{kreweras}, and
then rederived by Niederhausen~\cite{niederhausen1} and
Gessel~\cite{gessel-proba}. It is, however, not
well-understood, for two reasons:

-- no direct proof of~\Ref{kreweras-simple} is known, even when
$i=0$. The number of walks ending at the
origin is closely related to the number of non-separable
planar maps, to the number of {\em cubic\/} non-separable
maps~\cite{poulalhon-schaeffer,gilles-STOC,tutte1,tutte2}, and to 
the number of two-stack sortable
permutations~\cite{mbm-2ss,west-2ss,zeil-2ss}.  
All available proofs of~\Ref{kreweras-simple}  are rather long and 
complicated. Moreover, in all of them, the result is {\em checked\/} rather
than {\em derived\/}.

-- most importantly, the three-variate \gf \ for these walks can be
shown to be algebraic~\cite{gessel-proba}, but none of the proofs
explain combinatorially this algebraicity.

\medskip

All problems of walks confined in a quadrant can be attacked by
writing a functional equation for their three-variate \gf , and it
is this uniform approach that we discuss here. This
functional equation  simply encodes the step-by-step construction
of the walks. For instance, for square lattice walks, we can write
{ $$
Q(x,y;t):=\sum_{i,j,n \ge 0} a_{i,j}(n) x^i y^j t^n \hskip 8cm$$}
$$\hskip 16mm= 1 + t(x+y) Q(x,y;t) + t \ \frac{Q(x,y;t) - Q(0,y;t)} x
+ t \ \frac{Q(x,y;t) - Q(x,0;t)} y,
$$
that is, 
\begin{equation} 
 \left(xy -t(x+y+x^2y+xy^2)  \right) Q(x,y;t) = xy-xtQ(x,0;t)
-ytQ(0,y;t),
\label{qdp-carre} \end{equation}
and the solution of this equation, given by~\Ref{carre-general}, is
D-finite (but transcendental). Similarly, for the diagonal square
lattice, we have 
$$
\left( xy-t(1+x^2)(1+y^2)\right) Q(x,y;t) = \hskip 8cm$$
$$\hskip 3cm xy-t(1+x^2)Q(x,0;t)
-t(1+y^2)Q(0,y;t)+tQ(0,0;t), 
$$
with again a D-finite transcendental solution, while for Kreweras'
algebraic model, we obtain
\begin{equation}
\left( xy-t(x+y+x^2y^2)\right) Q(x,y;t) = xy-xtQ(x,0;t) -ytQ(0,y;t).
\label{qdp-kreweras}
\end{equation}
Finally,  the equation that governs  the non-holonomic model of~\cite{mbm-marko2} is
$$
\left( xy-t(x^3+y^3)\right) Q(x,y;t) = x^2y^2-tx^3Q(x,0;t)
-ty^3Q(0,y;t).
$$

The general theme of this paper is the following: the above equations
completely solve, in some sense, the problem of enumerating the
walks. But they are not the kind of solution one likes, especially if
the numbers are simple, or if the \gf \  is actually algebraic! How
can one derive these simple solutions from the functional equations?
And what is the essential difference between, say, Eqs.~\Ref{qdp-carre}
and~\Ref{qdp-kreweras}, that makes one series transcendental, and the
other algebraic?

We shall answer some of these questions. Our main result is a new
proof of~\Ref{kreweras-simple}, which we believe to be simpler than the
three previous ones. It has, at least, one nice feature: we {\em derive the
algebraicity from the equation\/} without having to  guess
the formula first.
Then, we give a new proof of a (refinement of) a holonomy criterion
that was proved combinatorially in~\cite{mbm-marko2}:
 if the set of steps $\frak S$ is symmetric with respect to the
$y$-axis and satisfies a {\em small horizontal variations\/} condition,
then the \gf \ for walks with steps in $\frak S$, starting from any
given point $(i_0,j_0)$, is D-finite. This result covers the two above
D-finite transcendental cases, but not Kreweras' model... We finally
survey some perspectives of this work.

\bigskip
Let us  conclude this section with a few more formal definitions
on walks and power series.

Let $\frak S$ be a finite subset of $\zs ^2$. A walk with steps in
$\frak S$ is a
finite sequence $w = (w_0, w_1, \ldots , w_n)$  of vertices of $\zs
^2$ such that $w_{i}-w_{i-1} \in \frak S$ for $1 \le i \le
n$. The number of 
steps, $n$, is the {\em length\/} of $w$. The  starting point
of $w$ is $w_0$, and its  endpoint is
$w_n$.  The {\em complete \gf\/} for  a set $\frak A$ of walks starting
from a given point $w_0 =(i_0,j_0)$ 
is the series
$$ A(x,y;t)=\sum_{n \ge 0}t^n \sum_{i,j \in \zss} a_{i,j}(n) x^i 
y^j ,$$ where $a_{i,j}(n)$ is the number  of walks of  $\frak A$
that have length $n$  and  end at $(i,j)$. This series is a formal
power series in $t$ whose coefficients are polynomials in $x, y, 1/x,
1/y$. We shall often denote $\bx =1/x$ and $\by =1/y$.

 Given a ring $\GL$ and $k$ variables
$x_1, \ldots , x_k$, we denote by
 $\GL[x_1, \ldots , x_k]$ the ring of polynomials in $x_1,
\ldots , x_k$ with coefficients in $\GL$,
and by $\GL[[x_1, \ldots , x_k]]$ the ring of  formal power
series in $x_1, 
\ldots , x_k$ with coefficients in $\GL$.
%
%
%
 If $\GL$ is a field, we denote by
 $\GL(x_1, \ldots , x_k)$ the field of rational functions in $x_1,
\ldots , x_k$ with coefficients in $\GL$.
%

%
%

Assume $\GL$ is a field. A series $F$ in $\GL[[x_1, \ldots , x_k]]$ is
{\em rational\/} if 
there exist  polynomials $P$ and $Q$
in $\GL[x_1, \ldots ,
x_k]$, with $Q\not = 0$, such that $QF=P$. 
It is {\em algebraic\/}
(over the field  $\GL(x_1, \ldots , x_k)$) if
 there exists a non-trivial polynomial $P$ with coefficients in
$\GL$ such that 
$P(F,x_1, \ldots , x_k)=0.$ The sum and product of algebraic series
is algebraic.  
%

The series $F$  is {\em D-finite\/} (or {\em holonomic\/})  
if the partial derivatives of $F$ span a finite 
dimensional vector space over the field $\GL(x_1, \ldots , x_k)$
(this vector space is a subspace  of the fraction field of $\GL[[x_1, \ldots , x_k]]$);
see~\cite{stanleyDF} for the one-variable case,
and~\cite{lipshitz-diag,lipshitz-df} otherwise.  In other 
words, for $1\le 
i\le k$, the series $F$ satisfies a non-trivial partial differential
equation of the form
$$\sum_{\ell=0}^{d_i}P_{\ell,i}\frac{\partial ^\ell
F}{\partial x_i^\ell} =0,$$
where $P_{\ell,i}$ is a polynomial in the $x_j$.
Any algebraic series is holonomic. The
sum and product of two 
holonomic series are still holonomic. The
specializations of a holonomic series (obtained by giving 
values from $\GL$  to some of the variables) are holonomic, if
well-defined. Moreover, if $F$ is an {\em algebraic\/} series and 
$G(t)$ is a holonomic series of one variable, then the substitution
$G(F)$ (if well-defined) is
holonomic~\cite[Prop.~2.3]{lipshitz-df}. 

\section{\bf A new proof of Kreweras' result}
Consider walks that start from $(0,0)$, are made of South, West and
North-East steps, and always stay in the first quadrant
(Figure~\ref{figure-kreweras}). Let 
$a_{i,j}(n)$ be the number of $n$-step walks of this type ending at $(i,j)$. We
denote by $Q(x,y;t)$ the complete \gf  \ of these walks:
$$ 
Q(x,y;t):=\sum_{i,j,n \ge 0} a_{i,j}(n) x^i y^j t^n .
$$
Constructing the walks step by step yields the following equation:
\begin{equation}
\left( xy-t(x+y+x^2y^2)\right) Q(x,y;t) = xy-xtQ(x,0;t) -ytQ(0,y;t).
\label{qdp-kreweras1}
\end{equation}
We shall often denote, for short, $Q(x,y;t)$ by $Q(x,y)$. Let us
also denote the series $xtQ(x,0;t)$ by $R(x;t)$ or even $R(x)$. Using the
symmetry of the problem in $x$ and $y$, the above equation
becomes:
\begin{equation}
\left( xy-t(x+y+x^2y^2)\right) Q(x,y) = xy-R(x) -R(y).
\label{qdp-kreweras2}
\end{equation}
This equation is equivalent to a recurrence relation defining the
numbers $a_{i,j}(n)$ by induction on $n$. Hence, it defines completely
the series $Q(x,y;t)$. Still, the characterization of this series we
have in mind is of a different nature:
\begin{Theorem}
\label{theorem-kreweras}
Let $X\equiv X(t)$ be the power series in $t$ defined by
$$X=t(2+X^3).$$
Then the \gf \ for Kreweras' walks ending on the $x$-axis is
$$Q(x,0;t)= \frac 1 {tx} \left( \frac 1 {2t} - \frac 1 x - 
\left( \frac 1 X -\frac 1 x \right) \sqrt{1-xX^2} \right).$$
Consequently, the length \gf \ for walks ending at $(i,0)$ is
$$[x^i] Q(x,0;t) = \frac{X^{2i+1}}{2.4^i\ t}\left( C_i
-\frac{C_{i+1}X^3}4\right),$$ 
where $C_i={{2i} \choose i}/(i+1)$ is the $i$-th Catalan number.
The Lagrange inversion formula gives the number of such walks of 
length $3n+2i$ as
$$a_{i,0}(3n+2i)=\frac{4 ^n (2i+1)}{(n+i+1)(2n+2i+1)} {{2i} \choose i}
{{3n+2i} \choose n} .$$
\end{Theorem}
The aim of this section is to derive Theorem~\ref{theorem-kreweras}
from the functional equation~\Ref{qdp-kreweras1}.

\medskip
\noindent{\bf Note.} Kreweras also gave a closed form expression for the number
of walks containing exactly $p$ West steps, $q$ South steps, and $r$
North-East steps, that is, for walks of length $m=p+q+r$ ending at
$(i,j)=(r-p,r-q)$:
$$a_{r-p,r-q}(p+q+r)= {{p+q+r} \choose{p,q,r}} \left( 1-
\frac{p+q}{r+1}\right)\hskip 8cm$$
$$ +\sum_{h=1}^p \sum_{k=1}^q 
\frac{(-1)^{h+k}}{(h+k)(h+k-1)}{{h+k}\choose h}
{{2h+2k-2}\choose{2h-1}} {{p+q+r} \choose{p-h,q-k,r+h+k}}.$$
The functional equation~\Ref{qdp-kreweras1}, combined with the expression of
$Q(x,0)$ given in Theorem~\ref{theorem-kreweras}, gives an alternative
expression for this number, 
still in the form of a double sum:
$$a_{r-p,r-q}(p+q+r)= {{p+q+r} \choose{p,q,r}}
- \sum_{i\ge 0}\sum_{n\ge 0}\frac{4 ^n (2i+1)}{(n+i+1)(2n+2i+1)} {{2i} \choose i}{{3n+2i} \choose n}$$
$$\hskip 1cm \times \left(  {{p+q+r-3n-2i-1}
\choose{p-n,q-n-i-1,r-n-i}}
+ {{p+q+r-3n-2i-1} \choose {p-n-i-1,q-n,r-n-i}}\right) .$$
This expression has a straightforward combinatorial explanation (all
walks, except those that cross the $x$- or $y$-axis). But none of these
formulas specialize to the above simple expression of 
$a_{i,0}(m)$ when $q=r$...
\subsection{The obstinate kernel method}
The kernel method is basically the only tool we have to attack
Equation~\Ref{qdp-kreweras2}. This method had been around
since, at least, the 
70's, and is currently the subject of a certain revival (see the
references
in~\cite{hexacephale,banderier-flajolet,bousquet-petkovsek}).
It consists in coupling the
variables $x$ and $y$ so as to cancel the {\em kernel\/}  $K(x,y)=
xy-t(x+y+x^2y^2)$.
This should give the ``missing'' information about
the series $R(x)$.

As a polynomial in $y$, this kernel has two roots
$$
\begin{array}{lclllll}
Y_0(x)&=&\displaystyle \frac{1-t\bx -\sqrt{(1-t\bx )^2-4t^2x}}{2tx}&=&
 & &t+\bx t^2 + O(t^3) , \\
\\
 Y_1(x)&=&\displaystyle \frac{1-t\bx +\sqrt{(1-t\bx )^2-4t^2x}}{2tx}
&=&\displaystyle\frac \bx t -\bx ^2 &-&t-\bx t^2 + O(t^3) .
\end{array}$$
The elementary symmetric functions of the $Y_i$ are
\begin{equation}
Y_0+Y_1= \frac \bx {t} - \bx ^2 \quad \hbox{and}
\quad Y_0Y_1 = \bx .
\label{symmetric}
\end{equation}
The fact that they are polynomials in $\bx =1/x$ will play a very
important role below.

Only the first root can be substituted for $y$
in~\Ref{qdp-kreweras2} (the term $Q(x,Y_1;t)$ is not a well-defined
power series in $t$). We thus obtain a functional equation for
$R(x)$:
\begin{equation}
R(x)+R(Y_0)=xY_0.
\label{kernel0}
\end{equation}
It can be shown that this equation -- once restated in terms of
$Q(x,0)$ -- defines uniquely $Q(x,0;t)$ as a
formal power series in $t$ with polynomial coefficients in
$x$. Equation~\Ref{kernel0} is the standard result of the kernel method.

Still, we want to apply here the {\em obstinate\/} kernel method. That
is, we shall not content ourselves with Eq.~\Ref{kernel0}, but we
shall go on producing pairs $(X,Y)$ that cancel the kernel and use the
information they provide on the series $R(x)$.
This obstination was inspired by the book~\cite{fayolle} by Fayolle,
Iasnogorodski and Malyshev, and more precisely by Section~2.4 of this
book, where one possible way to obtain such pairs is described (even
though the analytic context is different). We give here
an alternative construction that actually  provides the same pairs.

\begin{figure}[hbt]
\begin{center}
 \scalebox{1}{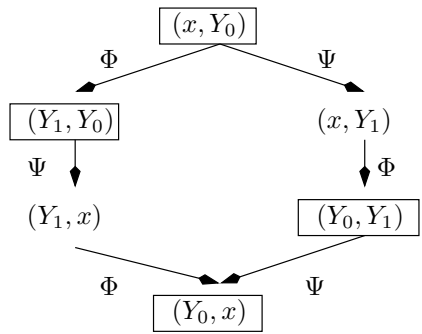 }
\end{center}
\caption{The orbit of  $(x,Y_0)$ under the action of  $\Phi$ and $\Psi$.}
\label{diagram}
\end{figure}

Let $(X,Y)\not = (0,0)$ be a pair of Laurent series in $t$ with
coefficients in a field $\GK$  such that $K(X,Y)=0$. 
We define $\Phi(X,Y)= (X',Y)$, where $X'=(XY)^{-1}$ is {\em the other
solution\/}  of $K(x,Y)=0$, seen as a polynomial  in $x$.
Similarly, we define $\Psi(X,Y)= (X,Y')$, where $Y'=(XY)^{-1}$ is 
the other solution  of $K(X,y)=0$. Note that $\Phi$ and $\Psi$ are
involutions. Let us examine their action on the
pair $(x,Y_0)$. We obtain the diagram of Figure~\ref{diagram}.

\label{page-diagram}
All these pairs of power series cancel the kernel, and we have framed
the ones that can be legally substituted\footnote{The fact that the
series $Q(Y_0,Y_1;t)$ and $Q(0,Y_1;t)$ {\em are\/} 
well-defined is not immediate, and depends strongly on the three steps
taken by the walks.} in the main functional 
equation~\Ref{qdp-kreweras2}. We thus obtain
 {\em two\/} equations for the unknown series $R(x)$:
\begin{eqnarray}
R(x)+R(Y_0)&=&xY_0, \label{kernel2.1}\\
R(Y_0)+R(Y_1)&=&Y_0Y_1=\bar x. \label{kernel2.2}
\end{eqnarray}

\medskip
\noindent{\bf Remark.} Let $p,q,r$ be three nonnegative numbers such that
$p+q+r=1$. Take $x= (pr)^{1/3}q^{-2/3}$,  $y= (qr)^{1/3}p^{-2/3}$, and
$t=(pqr)^{1/3}$. Then $K(x,y;t)=0$, so that $R(x)+R(y)=xy$. This
equation can be given a probabilistic interpretation by considering
random walks that make a North-East step with (small) probability $r$
and a South (resp. West) step  with probability $p$ (resp.~$q$). This
probabilistic argument, and the equation it implies, is the starting
point of Gessel's solution of Kreweras 
problem~\cite[Eq.~(21)]{gessel-proba}.

\subsection{Symmetric functions of $Y_0$ and $Y_1$}
After the kernel method, the next tool of our approach is the extraction
of the positive part of power series. More precisely, let $S(x;t)$ be
a power series in $t$ whose coefficients are Laurent polynomials in
$x$:
$$S(x;t)= \sum_{n\ge 0} t^n \sum_{i \in \zss} s_i(n) x^i ,$$
where for each $n \ge 0$, only finitely many coefficients $s_i(n)$ are
non-zero. We define the positive part of this series by
$$S^+(x;t):= \sum_{n\ge 0} t^n \sum_{i \in \nss} s_i(n) x^i .$$
This is where the values of the symmetric functions of $Y_0$ and $Y_1$
become crucial: the fact that they only involve negative powers of $x$
(see~\Ref{symmetric}) will simplify the extraction of the positive
part of certain equations.
\begin{Lemma}
\label{lemma-symmetric}
Let $F(u,v;t) $ be a power series in $t$ with coefficients in
$\cs[u,v]$, such that $F(u,v;t)=F(v,u;t)$. Then the series
$F(Y_0,Y_1;t)$, if well-defined, is a power series in $t$ with
polynomial coefficients in $\bx$. Moreover, the constant term of this
series, taken with respect to $\bx$, is $F(0,0;t)$.
\end{Lemma}
\noindent {\bf  Proof.} All symmetric polynomials of $u$ and $v$
are polynomials in $u+v$ and $uv$ with complex coefficients.
\cqfd

We now want to form a symmetric function of $Y_0$ and $Y_1$, starting
from the equations~(\ref{kernel2.1}--\ref{kernel2.2}). The first one
reads
$$R(Y_0)-xY_0=-R(x).$$
By combining both equations, we then obtain the companion expression:
$$R(Y_1)-xY_1=R(x)+ 2\bx -1/t.$$
Taking the product\footnote{An alternative derivation of Kreweras'
result, obtained by considering the divided difference
$(R(Y_0)-xY_0-R(Y_1)+xY_1)/(Y_0-Y_1)$, will be discussed on the
complete version of this paper.} of these two equations gives
$$(R(Y_0)-xY_0)(R(Y_1)-xY_1)=-R(x)(R(x)+ 2\bx -1/t).$$
The extraction of  the positive part of this identity is made possible
by Lemma~\ref{lemma-symmetric}. Given that $R(x;t)=xtQ(x,0;t)$, one obtains:
$$x= -t^2x^2Q(x,0)^2 +(x-2t)Q(x,0)+2tQ(0,0),$$
that is,
\begin{equation}
t^2x^2Q(x,0)^2 +(2t-x)Q(x,0)-2tQ(0,0)+x=0.
\label{quadratic}
\end{equation}

\subsection{The quadratic method}

Equation~\Ref{quadratic} -- which begs for a combinatorial explanation
-- is typical of the equations obtained when
enumerating planar maps, and the rest of the proof will be routine to all maps
lovers. This equation can be solved using the so-called {\em quadratic
method\/}, which
was first invented by Brown~\cite{brown}. The formulation we use here
is different 
both from Brown's original presentation and from the one in Goulden and
Jackson's book~\cite{gj}. This new formulation is convenient for
generalizing the method to equations of higher degree with more
unknowns~\cite{mbm-phoenix}. 

Equation~\Ref{quadratic} can be written as 
\begin{equation} P(Q(x), Q(0),t,x)=0,\label{quad1}\end{equation}
where 
$P(u,v,t,x)= t^2x^2u^2 +(2t-x)u-2tv+x$,
and $Q(x,0)$ has been abbreviated in $Q(x)$. Differentiating this
equation with respect to $x$, we find
$$
\frac{\partial P}{\partial u}(Q(x), Q(0),t,x)\frac{\partial
Q}{\partial x}(x)
+ \frac{\partial P}{\partial x}(Q(x), Q(0),t,x)=0.
$$
Hence,  if there exists a power series in $t$, denoted $X(t)\equiv X$, such that
\begin{equation}
\frac{\partial P}{\partial u}(Q(X), Q(0),t,X)=0,
\label{quad2} \end{equation}
then one also has
\begin{equation}  
 \frac{\partial P}{\partial x}(Q(X), Q(0),t,X)=0,
\label{quad3} \end{equation}          
and we thus obtain a system of three polynomial equations, namely
Eq.~\Ref{quad1} 
written for $x=X$, Eqs.~\Ref{quad2} and~\Ref{quad3}, that relate
the three unknown series $Q(X)$, $Q(0)$ and $X$. This puts us in a
good position to write an algebraic equation defining $Q(0)=Q(0,0;t)$.

Let us now work out the details of this program: Eq.~\Ref{quad2} reads
$$X=2t^2X^2Q(X) +2t,$$
and since the right-hand side is a multiple of $t$, it should be clear
that this equation defines a unique power 
series $X(t)$. The system of three equations reads
$$
\left\{
\begin{array}{l}
t^2X^2Q(X)^2 +(2t-X)Q(X)-2tQ(0)+X=0,\\
2t^2X^2Q(X) +2t-X=0,\\
2t^2XQ(X)^2 -Q(X)+1=0.
\end{array}
\right.
$$
Eliminating $Q(X)$ between the last two equations yields
$X=t(2+X^3),$
%
so that the series $X$ is the parameter introduced in
Theorem~\ref{theorem-kreweras}.  
Going on with the elimination, we finally obtain
$$Q(0,0;t)=\frac X{2t} \left( 1-\frac {X^3} 4\right),$$
and the expression of $Q(x,0;t)$ follows from~\Ref{quadratic}.
\cqfd

\section{\bf A holonomy criterion}
Using functional equations, we can recover, and actually refine, a
holonomy criterion that was recently proved
combinatorially~\cite{mbm-marko2}. 
Let $\frak S$ be a finite subset of $\zs^2$. We say that  $\frak S$ is
symmetric with respect to the $y$-axis if 
$$(i,j) \in {\frak S} \Rightarrow (-i,j) \in {\frak S}.$$
We say that $\frak S$ has small horizontal variations if 
$$(i,j) \in {\frak S} \Rightarrow |i|\le 1.$$
The usual square lattice steps satisfy these two conditions. So do the
steps of the diagonal square lattice (Figure~\ref{diagonal}).

\begin{Theorem}\label{theorem-sufficient}
Let ${\frak S}$   be a finite subset of $\zs^2$ that is symmetric with
respect to the $y$-axis and has small horizontal variations. Let
$(i_0,j_0) \in \ns^2$. Then the complete \gf \ $Q(x,y;t)$ for walks
that start from 
$(i_0, j_0)$, take their steps in $\frak S$ and stay  in the first
quadrant is D-finite.
\end{Theorem}
A combinatorial argument proving the holonomy of $Q(1,1;t)$ is presented
in~\cite{mbm-marko2}. 

\subsection{Example}
Before we embark on the proof of this theorem, let us see the
principle of the proof
at work  on a simple example: square lattice walks confined in a
quadrant.
The functional equation satisfied by their complete \gf \ is 
$$ \left(xy -t(x+y+x^2y+xy^2)  \right) Q(x,y) = \hskip 8cm$$
\vskip -6mm
\begin{equation} 
\hskip 40mm 
xy-xtQ(x,0)
-ytQ(0,y) =xy-R(x)-R(y),
\label{eq-carre}
\end{equation}
where, as in Kreweras' example, we denote by $R(x)$ the series
$txQ(x,0)$. The kernel $K(x,y)=xy 
-t(x+y+x^2y+xy^2)$, considered as a polynomial in $y$, has two roots:
$$Y_0(x)=\displaystyle \frac{1-t(x+\bx) -\sqrt{(1-t(x+\bx) )^2-4t^2}}{2t}=
 \hskip 18mm t+(x+\bx) t^2 + O(t^3) ,$$
 $$ Y_1(x)=\displaystyle \frac{1-t(x+\bx) +\sqrt{(1-t(x+\bx) )^2-4t^2}}{2t}
=\displaystyle\frac 1 t -x-\bx  -t-(x+\bx) t^2 + O(t^3) .$$
The elementary symmetric functions of the $Y_i$ are
$$
Y_0+Y_1= \frac 1 {t} - x-\bx \quad \hbox{and}
\quad Y_0Y_1 = 1 .
$$
Observe that they are no longer polynomials in $\bx =1/x$.

If, as above, we apply to the pair $(x,Y_0)$ the transformations $\Phi$ and
$\Psi$, we obtain a very simple diagram:

\begin{figure}[hbt]
\begin{center}
 \scalebox{1}{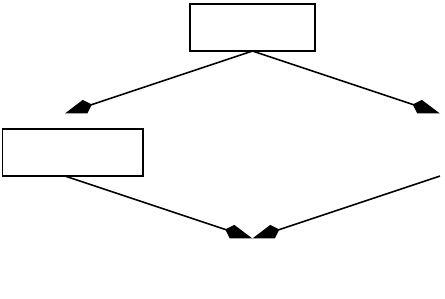 }
\end{center}
\end{figure}

\noindent 
From the two pairs that can be substituted for $(x,y)$ in Equation~\Ref{eq-carre},
we derive the following system:
\begin{eqnarray*}
R(x)+R(Y_0)&=&xY_0, \\
R(\bx)+R(Y_0)&=&\bar x Y_0. 
\end{eqnarray*}
From here, the method has to diverge from what we did in Kreweras'
case. Eliminating $R(Y_0)$ between the two equations gives
\begin{equation}
R(x)-R(\bx) = (x-\bx ) Y_0.
\label{diff-sym}
\end{equation}
Since $R(0)=0$,  extracting the positive part of this
identity gives $R(x)$ as {\em the positive part of an algebraic
series\/}. It is known that the positive part of a D-finite
series is always D-finite~\cite{lipshitz-diag}. In particular, the
series $R(x)$ is D-finite. The same holds for $Q(x,0)$, and,
by~\Ref{eq-carre}, for $Q(x,y)$. 

This argument is enough for proving the holonomy of the series, but,
given the simplicity of this model, we can proceed with explicit
calculations. Given the polynomial equation defining $Y_0$,
$$Y_0=t(1+\bx Y_0 + xY_0+Y_0^2)=t(1+\bx Y_0)(1+x Y_0),$$
the Lagrange inversion formula yields the following expression for
$Y_0$:
$$Y_0= \sum_{m\ge 0} \sum_{i \in \zss} \frac{x^i\
t^{2m+|i|+1}}{2m+|i|+1}
{{2m+|i|+1} \choose {m+|i|}}
{{2m+|i|+1} \choose {m}}.$$
Since $R(0)=0$,  extracting the positive part in the
identity~\Ref{diff-sym} now gives, after some reductions,
$$R(x)=txQ(x,0) = 
\sum_{m\ge 0} \sum_{i \ge 0} \frac{x^{i+1}
t^{2m+i+1}(i+1)}{(2m+i+1)(2m+i+2)}
{{2m+i+2} \choose {m+i+1}}
{{2m+i+2} \choose {m}}.$$
This naturally fits with the general expression~\Ref{carre-general}.

\subsection{Proof of Theorem~\ref{theorem-sufficient}.}
\noindent 
We define two Laurent polynomials in $y$ by
$$P_0(y) := \sum_{(0,j) \in {\frak S}} y^j \quad \hbox{and}
\quad P_1(y) := \sum_{(1,j) \in {\frak S}}y^j.$$ Let $-p$ be the
largest down move; more precisely, 
$$p=\max(0,\{-j : (i,j) \in {\frak S} \ \hbox{for some } i\}).$$ 
 The functional equation obtained by constructing walks step-by-step reads:
$$
 K( x,y) Q( x,y) = \hskip 12cm
$$
\vskip -6mm
\begin{equation} x^{1+i_0}y^{p+j_0}
 -ty^pP_1(y) Q(0,y)
- t\sum_{(i,-j) \in {\frak S}}  \ 
\sum_{m=0}^{j-1} \left( Q_m(x) - \delta_{i,1} Q_m(0)\right) x^{1-i} y^{p+m-j}
 \label{func-eq} 
\end{equation}
where
$$K(x ,y)= xy^p\left( 1- tP_0(y)-t(x+\bx) P_1(y)\right) $$ is the
 kernel of the equation, and $Q_m(x)$ stand for the coefficient
of $y^m$ in $Q(x,y)$. All the series involved in this equation also
depend on the variable $t$, but it is omitted for the sake of briefness.
For instance, $K(x,y)$ stands for $K(x,y;t)$.

 As above,
we shall use the kernel method -- plus another argument -- to solve
the above functional equation. The polynomial $K(x,y)$, seen as a
polynomial in $y$, admits a number of roots, which are Puiseux series
in $t$ with coefficients in an algebraic closure of $\qs(x)$.
 Moreover, all these roots are distinct.
 As $K(x,y;0)=xy^p$, exactly $p$ of these roots, say $Y_1,
\ldots , Y_p$, vanish at $t=0$. This property guarantees that these
$p$ series can be substituted for $y$ in~\Ref{func-eq}, which yields
\begin{equation} 
 x^{1+i_0}Y^{p+j_0}=
 tY^pP_1(Y) Q(0,Y)
+ t\sum_{(i,-j) \in {\frak S}}  \ 
\sum_{m=0}^{j-1} \left( Q_m(x) - \delta_{i,1} Q_m(0)\right) x^{1-i} \ Y^{p+m-j}
 \label{kernel1}\end{equation}
for any $Y=Y_1, \ldots , Y_p$. 

Given the symmetry of $K$ in $x$ and $\bx$, each of the $Y_i$ is
invariant by the transformation $x \rightarrow 1/x$. Replacing $x$ by
$\bar x$ in the above equation gives, for any $Y=Y_1, \ldots , Y_p$,
\begin{equation} 
 \bx^{1+i_0}Y^{p+j_0}=
 tY^pP_1(Y) Q(0,Y)
+ t\sum_{(i,-j) \in {\frak S}}  \ 
\sum_{m=0}^{j-1} \left( Q_m(\bx) - \delta_{i,1} Q_m(0)\right) x^{i-1} Y^{p+m-j}
 .\label{kernel2}\end{equation}
We now combine \Ref{kernel1} and \Ref{kernel2} to eliminate $Q(0,Y)$:
$$ 
(x^{1+i_0}-\bx^{1+i_0})Y^{p+j_0}=
 t\sum_{(i,-j) \in {\frak S}}  \ 
\sum_{m=0}^{j-1} \left( x^{1-i}Q_m(x) - x^{i-1} Q_m(\bx)\right) Y^{p+m-j}
$$ 
for any $Y=Y_1, \ldots , Y_p$. This is the generalization of
Eq.~\Ref{diff-sym}. The right-hand side of the above
equation is a polynomial $P$ in $Y$, of degree at most $p-1$. We know
 its value at 
$p$ points, namely $Y_1, \ldots , Y_p$. The Lagrange
interpolation formula implies that these $p$ values completely
determine the polynomial. As the left-hand side of the equation is 
 algebraic, then 
each of the coefficients of $P$ is also algebraic. That is,
$$ 
 t\sum_{(i,-j) \in {\frak S}}  \ 
\sum_{m=0}^{j-1} \left( x^{1-i}Q_m(x) - x^{i-1} Q_m(\bx)\right) y^{p+m-j}
= \sum_{m=0}^{p-1} A_m(x) y^m,$$
where each of the $A_m$ is an algebraic series.
 Let us extract the positive part of this
 identity. Given that $i$ can only be $0$, $1$ or $-1$, we obtain
$$   t\sum_{(i,-j) \in {\frak S}}  \ 
\sum_{m=0}^{j-1} \left( x^{1-i}Q_m(x) -  \delta_{i,1}
Q_m(0)\right) y^{p+m-j}
= \sum_{m=0}^{p-1} H_m(x) y^m
$$
where $ H_m(x):=A_m^+(x)$ is the positive part of $A_m(x)$. Again,
this series can be shown to be D-finite.
Going back to the original functional equation~\Ref{func-eq}, this gives
$$ K(x,y) Q(x,y)
=  x^{1+i_0}y^{p+j_0}
  -ty^pP_1(y) Q(0,y)
-\sum_{m=0}^{p-1} H_m(x) y^m
.$$
Let us finally\footnote{In the square lattice case, the symmetry of the model in $x$ and
$y$ makes this step unnecessary: once the holonomy of $Q(x,0)$ is
proved, the holonomy of $Q(x,y)$ follows.} consider the kernel $K(x,y)$ as a  polynomial in
$x$. One of its roots, denoted
below $X$, is a formal power series in 
$t$ that vanishes at $t=0$. Replacing $x$ by this root allows us to
express $Q(0,y)$ as a D-finite series:
$$ ty^pP_1(y)Q(0,y) = X^{1+i_0}y^{p+j_0}
-\sum_{m=0}^{p-1} H_m(X) y^m.$$ 
 The functional
equation finally reads
$$K(x,y) Q(x,y)
=  (x^{1+i_0}-X^{1+i_0})y^{p+j_0}
 -\sum_{m=0}^{p-1} (H_m(x)- H_m(X)) y^m.$$
Since the substitution of an algebraic series into a D-finite
one gives another D-finite series, this equation shows that $Q(x,y)$
is D-finite. 
\cqfd

\section{\bf Further comments, and  perspectives}

We first give some asymptotic estimates for the number of $n$-step
walks in the quadrant, for various sets of steps. Then, a number of
research directions, 
which I  have started to explore, or would like to explore in 
the coming  months, are presented. All of them are motivated by the
new proof of Kreweras' formula given in Section~2.

\subsection{Asymptotic estimates}
Following a suggestion of one of the referees, the table below
summarizes the asymptotic behaviour on the number of $n$-step walks in
the first quadrant, for the four models mentioned in the introduction, with
three different conditions: the endpoint is fixed, the endpoint lies
on the $x$-axis, the endpoint is free. This is all the more relevant
that the argument proving the transcendence of the square lattice case
is based on asymptotic estimates. The results for the two versions of
the square lattice
can be obtained directly from the formulas given in the
introduction. The results for Kreweras' walks are derived from
Theorem~\ref{theorem-kreweras} and the functional
equation~\Ref{qdp-kreweras1}, by analysing 
the singularities of the series~\cite{flajolet-odlyzko}. The last series of
results is derived from~\cite{mbm-marko2} using, again, an analysis of the
singularities of the \gfs .

\bigskip
\hskip -5mm
\begin{tabular}{c|c|c|c|c}
Model	& Specific	& Endpoint	& 
Free & Nature \\
& endpoint	&  on the $x$-axis & endpoint	&of the series \\ \hline
&&&&\\
Ordinary or 	& $\displaystyle \frac{4^n}{n^3}$
&$\displaystyle \frac{4^n}{n^2}$&$\displaystyle \frac{4^n}{n}	$&	D-finite\\
diagonal  sq. lattice &&&& transcendental\\
	&&&&\\

Kreweras' walks 	&$ \displaystyle \frac{3^n}{n^{5/2} }$ 	&$\displaystyle \frac{3^n}{n^{7/4} }	$&$\displaystyle \frac{3^n}{n^{3/4}}$&algebraic	\\
&&&&\\
Knight walks 	&$0$ 	&$\displaystyle
\frac{1}{n^{3/2}}\left(\frac{3}{4^{1/3}}\right)^n 	$&$2^n$	&
not D-finite\\

\end{tabular}

\subsection{Other starting points} It was observed by Gessel
in~\cite{gessel-proba} that the method he used to prove Kreweras'
result was hard to implement for a starting point different from the
origin. The reason of this difficulty is that, unlike the method
presented here, Gessel's approach {\em checks\/} the known expression
of the \gf , but does not {\em construct\/} it. I am confident that
the new approach of Section~2 can be used to solve such questions. If
the starting point does not lie on the main diagonal, the $x-y$
symmetry is lost; the diagram of Figure~\ref{diagram} now  gives
{\em four\/} different equations between the 
{\em two\/} unknown functions $Q(x,0)$ and $Q(0,y)$. 

\subsection{Other algebraic walk models}
\begin{figure}[ht]
\begin{center}
\includegraphics[height=30mm]{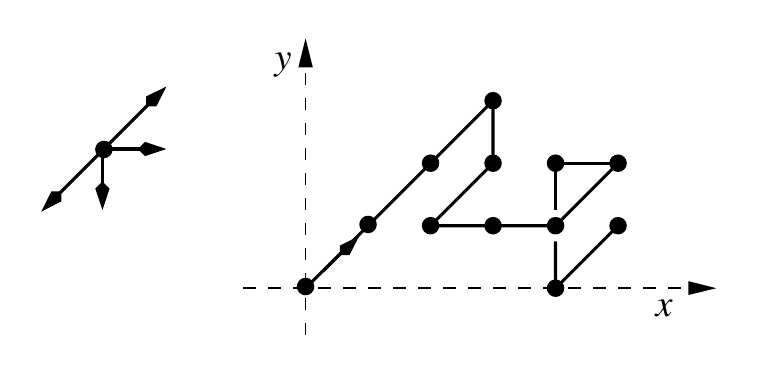}
\end{center}
\caption{A degenerate set of steps.}
\label{figure-degenerate}
\end{figure}

A close examination of the ingredients that make the proof of
Section~2 work might help to construct other walk models which, for
non-obvious reasons, would have an algebraic \gf . Note that for some
degenerate sets of steps, like those of
Figure~\ref{figure-degenerate}, the quadrant condition is equivalent
to a half-plane condition and thus yields an algebraic series. 

I have started a systematic exploration of walks with few steps and only one
up step: the non-trivial algebraic cases do not seem to be legion!
However, I met in this exploration one model that seems to yield nice
numbers (with a D-finite  \gf ) and for which the method
of Section~2 ``almost'' works. I then realized that the same problem
had been communicated to me, under a slightly different form, by Ira
Gessel, a few months ago. I plan to explore this model further.

\subsection{Other equations} Any combinatorial problem that seems to
have an algebraic    \gf \ and for which a linear functional equation
with two ``catalytic'' variables (in the terminology of
Zeilberger~\cite{zeil-umbral}) is available is now likely to be
attacked via the 
method of Section~2. These conditions might seem very restrictive, but
there is at least one such problem! The {\em vexillary\/} involutions,
which were conjectured in 1995 to be counted by Motzkin numbers,
satisfy the following equation:
$$\left( 1+\frac{t^2x^2y}{1-x} +\frac{t^2y}{1-y}\right)
F(x,y;t)=\hskip 8cm$$
$$\hskip 30mm \frac{t^2x^2y^2}{(1-ty)(1-txy)}+t\left(1+\frac{ty}{1-y}\right)F(xy,1;t)
+\frac{t^2x^2y}{1-x} F(1,y;t).$$
The conjecture was recently proved via a difficult combinatorial
construction~\cite{guibert-pergola-pinzani}. I have been able to apply
successfully the method of Section~2 to this
equation~\cite{mbm-vexillaires}.

\subsection{Random walks in the quarter plane}
Random walks in the quarter plane are naturally studied in
probability. Given a Markov chain on the first quadrant, a central
question is the determination of an/the invariant measure  $(p_{i,j})_{i,j
\ge 0}$. The invariance is equivalent to a linear equation satisfied
by  the series $P(x,y)
= \sum p_{i,j} x^i y^j$, in which the variables
$x$ and $y$ are ``catalytic''. A whole recent book is devoted to the
solution of this equation in the case where the walk has small
horizontal and vertical variations~\cite{fayolle}.
This book contains {\em one\/} example for which the series $P(x,y)$
is algebraic: no surprise, the steps of the corresponding walk are
exactly Kreweras' steps... This result is actually due to
Flatto and  Hahn~\cite{flatto}. The equation satisfied by the series
$P(x,y)$ does not work exactly like the equations for complete \gfs \
like $Q(x,y;t)$: roughly speaking, the third variable $t$ is replaced
by the additional constraint  $P(x,y)=1$.

Very recently, I have found  a new, simpler proof
of Flatto and Hahn's result (at least, in the symmetric case). 
The principle is the same as in
Section~2. One can either study a  
version of the enumeration problem in which each walk is weighted by
its probability (so that the invariant
distribution is obtained as a limit distribution), or directly adapt
the method  to the context of the series $P(x,y)$.
With both approaches, one remains, from the beginning to the end of
the proof, in the field of
algebraic series~\cite{mbm-algebraic}. This offer a significant
shortcut to Flatto and 
Hahn's proof, which is based on non-trivial
complex analysis, and uses a parametrisation of the roots of the 
kernel by elliptic
functions, which are {\em not\/} algebraic.

\bigskip
\noindent
{\bf Acknowledgements.} To my shame, I must recall that, in the
lecture that I gave at FPSAC'01 in Phoenix, I mentioned  (part
of) Kreweras' result as a conjecture. I am very grateful to Ira Gessel
who enlightened my ignorance by giving me the right references.
I also thank the anonymous referees for their very valuable comments.




\end{document}